\documentclass[11pt]{article}
\usepackage[utf8]{inputenc}
\usepackage[T1]{fontenc}
\usepackage{lmodern}
\usepackage[a4paper]{geometry}
\usepackage[english]{babel}
\usepackage{pifont}
\usepackage{extsizes}
\usepackage{xcolor}
\usepackage{fancybox}
\usepackage{graphicx}
\usepackage{sectsty}
\usepackage{enumerate}
\usepackage{multicol}
\usepackage[english]{varioref}
\usepackage[fleqn]{amsmath}
\usepackage{mathtools}
\usepackage{amssymb}
\usepackage{amsmath}
\usepackage{mathrsfs}
\usepackage{bm}
\usepackage{amsthm}
\usepackage{blkarray}
\usepackage{subfigure}
\usepackage{lipsum}
\usepackage{stmaryrd}
\usepackage[babel=true]{csquotes}
\theoremstyle{break}
\newtheorem{Th}{Theorem}
\newtheorem{Prop}[Th]{Proposition}

\newcommand{\conv}[2][n]{\underset{#1\rightarrow #2}{\longrightarrow}}
\newcommand{\eq}[2][n]{\underset{#1\rightarrow #2}{\sim}}

\newtheorem{Le}{Lemma}

\newcommand{\PP}{\mathbb{P}}

\newcommand{\RR}{\mathbb{R}}
\newcommand{\ZZ}{\mathbb{Z}}
\newcommand{\NN}{\mathbb{N}}

\newcommand{\ind}[1]{\mathbf{1}_{#1}\,}

\newcommand{\EEE}[1]{\mathbb{E}\left[#1 \right]}
\newcommand{\PPP}[1]{\mathbb{P}\left(#1 \right)}

\newenvironment{prooft}[1]{\vskip 2mm\noindent {\bf Proof of #1.}}
                    {\hfill $\square$ \vskip 2mm \noindent}
\date{}  
\title{Some properties on extremes for transient random walks in random sceneries}
\author{Nicolas Chenavier\footnote{Universit\'e du Littoral C\^ote d'Opale, Laboratoire de Math\'ematiques Pures et Appliqu\'ees J. Liouville, France. Mail: nicolas.chenavier@univ-littoral.fr}, Ahmad Darwiche \footnote{Universit\'e du Littoral C\^ote d'Opale, Laboratoire de Math\'ematiques Pures et Appliqu\'ees J. Liouville, France. Mail: darwich.ahmad.92@gmail.com}}

\begin{document}
\maketitle
\begin{abstract}
Let $(S_n)_{n \geq 0}$ be a transient random walk in the domain of attraction of a stable law and let $(\xi(s))_{s \in \ZZ}$ be a stationary sequence of random variables. In a previous work, under conditions of type $D(u_n)$ and $D'(u_n)$, we established a limit theorem for the maximum of the first $n$ terms of the sequence $(\xi(S_n))_{n\geq 0}$ as $n$ goes to infinity. In this paper we show that, under the same conditions and under a suitable scaling, the point process of exceedances converges to a Poisson point process. We also give some properties of $(\xi(S_n))_{n\geq 0}$.  \\

\noindent\textbf{Keywords:} extreme values, random walks, point processes.\\ 

\noindent\textbf{Mathematics Subject Classification:} 60G70, 60G50, 60G55.
\end{abstract}

\section{Introduction}
In 2009, Franke and Saigo \cite{franke_saigo2009bis, franke_saigo_2009} considered the following problem. Let $(X_k)_{k\geq 1}$ be a sequence of centered, integer-valued i.i.d. random variables and let $S_0=0$ a.s. and $S_{n}=X_{1}+\dots+X_{n}$, $n\geq 1$. Assume that, for any $x \in \RR$, 
\begin{align*}
    \PPP{\frac{S_n}{n^{1/\alpha}} \leq x} \conv[n]{\infty} F_\alpha(x),
\end{align*}
where $F_\alpha$ is the distribution function of a stable law with characteristic function given by
\begin{align*}
    \phi(\theta)=\exp(-|\theta|^\alpha(C_1+iC_2 \text{sgn}\,\theta)), \quad \alpha \in (0,2].
\end{align*}
Let $(\xi(s))_{s\in \ZZ}$ be a stationary sequence of $\RR$-valued random variables which are independent of the sequence  $(X_k)_{k\geq 1}$. The sequence $(\xi(S_n))_{n\geq 0}$ is referred to as a \textit{random walk in a random scenery}. In \cite{franke_saigo_2009}, Franke and Saigo derive limit theorems for the random variable $\max_{i\leq n}\xi(S_i)$ as $n$ goes to infinity when the $\xi(s)$'s are i.i.d.. The statements of their theorems depend on the value of $\alpha$. When $\alpha< 1$ (resp. $\alpha> 1$), it is known that the random walk $(S_n)_{n\geq 0}$  is transient (resp. recurrent) \cite{KS,LeGRo}. An important concept concerning random walks is the \textit{range}. The latter is defined as the number of sites visited by the first $n$ terms of the random walk, namely $R_n:=\#\{S_1, \dots, S_n\}$. The following result, due to Le Gall and Rosen \cite{LeGRo}, deals with its asymptotic behavior.
\begin{Th}[LeGall and Rosen] 
\label{th:legallrosen}
\begin{enumerate}[(i)]
    \item If $\alpha<1$, then
\[
 \frac{R_{[nt]}}{n} \conv[n]{\infty}qt \quad \PP-a.s.  
\]
    with $q:=\PPP{S_k\neq 0, \forall k \geq 1}$.
    \item If $\alpha=1$, then
\[
    \frac{h(n)R_{[nt]}}{n} \conv[n]{\infty}t \quad \text{in} \quad L^p(\PP),
\]
where $h(n):=1+\sum_{k=1}^n \PPP{S_k=0}$.
\item If $1<\alpha\leq 2$, then for any $L\in \NN$ and any $t_1<\dots<t_L$, 
\[
 \frac{1}{n^{1/\alpha}}\left(R_{\lfloor nt_1\rfloor},\dots,R_{\lfloor nt_L\rfloor}\right) \conv[n]{\infty} \left(m(Y(0,t_1)),\dots,m(Y(0,t_L))\right),  
\]
in distribution.
\end{enumerate}
\end{Th}
In the above result, $\{Y(t), t \in \RR\}$ denotes the right-continuous $\alpha-$stable Lévy process with characteristic function given by $\phi(t \theta)$ and  $m$ is the Lebesgue measure on $\RR$. One of the results of \cite{franke_saigo_2009} is the following. If $u_n$ is a threshold such that $n\PPP{\xi>u_n}\conv[n]{\infty}\tau$ for some $\tau>0$, with $\xi=\xi(1)$, and if the $\xi(s)$'s are i.i.d. then \[\PPP{\max_{i\leq n}\xi(S_i)\leq u_n} \conv[n]{\infty} e^{-\tau q}\] for $\alpha <1$. Such a result was generalized in \cite{Nicolas_Ahmad} for sequences $(\xi(s))_{s\in \ZZ}$ which are not necessarily i.i.d., but which satisfy a slight modification of the classical $D(u_n)$ and $D'(u_n)$ conditions of Leadbetter (see \cite{L1,L2} for a statement of these conditions). 

In this paper, we give a more precise treatment of the extremes of $(\xi(S_n))_{n\geq 0}$. To do it, we assume that the threshold is of the form $u_n=u_n(x)=a_nx+b_n $ ($a_n \in \RR$, $b_n >0$
and $x\in \RR$) and that, for any $x\in \RR$, the following term exists and is finite: 
\begin{equation}
\label{eq:defnu}
\nu(x,\infty):=\lim_{n \to \infty} n\PPP{\xi>u_n(x)}.   
\end{equation}
The quantity $\nu$ defines a measure on some topological space $E$.  According to the Gnedenko's theorem  \cite{Gnedenko}, if $\xi$ is in the domain of attraction of an extreme value distribution $G$, then $\nu$ is of the form:
\begin{equation*}
\nu(x,\infty)=\begin{cases}x^{-\beta}, &    E=(0, \infty] \quad \text{if  $G$ is a Fréchet distribution};\\
(-x)^{-\delta}, &  E=(-\infty,0]\quad \text{if  $G$ is a Weibull distribution};  \\
e^{-x}, &  E=(-\infty, \infty]\quad \text{if  $G$ is a Gumbel distribution};
\end{cases}
\end{equation*}
for some $\beta, \delta>0$. Notice that if $P_n$ denotes the distribution of $\frac{\xi - a_n}{b_n}$, then \eqref{eq:defnu} can be rephrased as 
\begin{equation} 
\label{eq:defnu2}
nP_n(A) \conv[n]{\infty} \nu(A),
\end{equation}
for any Borel subset $A\subset \RR$. 
Secondly, we assume that the (stationary) sequence $(\xi(s))_{s\in \ZZ}$ satisfies conditions of type $D(u_n)$ and $D'(u_n)$ in the same spirit as in \cite{Nicolas_Ahmad}. To introduce the first one, we write  for each $i_1< \cdots < i_p$ and for each $u\in \RR$,  \[F_{i_{1},\ldots,i_{p}}(u) = \PPP{\xi(i_1)\leq u, \ldots, \xi(i_p)\leq u}.\]  

\paragraph{$\mathbf{D}(u_n)$ condition}  We say that $(\xi(s))_{s\in \ZZ}$ satisfies the $\mathbf{D}(u_{n})$ condition if there exist a sequence $(\alpha_{n,\ell})_{(n,\ell)\in \NN^2}$ and a sequence $(\ell_n)$ of positive integers such that $\alpha_{n,\ell_n}\conv[n]{\infty} 0$, $\ell_n=o(n)$, and 
\begin{equation*} 
|F_{i_{1},\ldots,i_{p}, j_{1},\ldots, j_{p'}}(u_{n})-F_{i_{1},\ldots,i_{p}}(u_{n})F_{j_{1},\ldots, j_{p'}}(u_{n})|\leq \alpha_{n,\ell}
\end{equation*}
 for any integers $i_{1}<\dots <i_{p}<j_{1}< \dots<j_{p'}$ such that $j_{1}-i_{p}\geq \ell$. Notice that the bound holds uniformly in $p$ and $p'$. Roughly, the  $\mathbf{D}(u_n)$ condition (see e.g. p29 in \cite{lucarini})  is a weak mixing property for the tails of the joint distributions.

The $\mathbf{D'}(u_{n})$ condition  (see e.g. p29 in \cite{lucarini}) is a local type property and precludes the existence of clusters of exceedances. To introduce it, we consider a sequence $(k_n)$ such that
\begin{equation}
\label{eq:defrn}
k_n\conv[n]{\infty}\infty, \quad  \frac{n^2}{k_n} \alpha_{n,\ell_n}\conv[n]{\infty}0, \quad k_n\ell_n=o(n),
\end{equation}
where $(\ell_n)$ and $(\alpha_{n,l})_{(n,l)\in \NN^2}$ are the same as in the $\mathbf{D}(u_{n})$ condition.

\paragraph{$\mathbf{D'}(u_{n})$ condition}
In conjunction with the $\mathbf{D}(u_n)$ condition, we say that $(\xi(s))_{s\in \ZZ}$ satisfies the $\mathbf{D'}(u_{n})$ condition  if there exists a sequence of integers $(k_n)$ satisfying \eqref{eq:defrn} such that 
\[\lim\limits_{n\rightarrow \infty} n \sum_{s=1}^{\lfloor n/k_n \rfloor}\PPP{\xi(0)>u_{n},\xi(s)>u_n}=0.\]
In the classical literature, the sequences $(\alpha_{n,l})_{(n,l)\in \NN^2}$ and $(k_n)$ only satisfy $k_n\alpha_{n,\ell_n}\conv[n]{\infty}0$ (see e.g. (3.2.1) in \cite{lucarini}) whereas in \eqref{eq:defrn} we have assumed that  $\frac{n^2}{k_n} \alpha_{n,\ell_n}\conv[n]{\infty}0$. In this sense, the $\mathbf{D'}(u_{n})$ condition as written above is slightly more restrictive than the usual $D'(u_n)$ condition. \\

Our paper is organized as follows. In Section \ref{sec:ppexceedances}, we prove that under suitable scaling the so-called point process of exceedances converges to a Poisson point process in the transient case. In Section \ref{sec:properties}, we give some properties of the random walk in random scenery. More precisely, we show that the (stationary) sequence $(\xi(S_n))_{n\geq 0}$ satisfies the classical $D(u_n)$ condition of Leadbetter, but does not satisfy the $D'(u_n)$ condition. Our results generalize \cite{franke_saigo_2009} for sequences $(\xi(s))_{s\in \ZZ}$ which are not i.i.d. but which only satisfy the  $\mathbf{D}(u_{n})$ and  $\mathbf{D'}(u_{n})$ conditions. We also give some  remarks on the so-called extremal index and on the $D^{(k)}(u_n)$ condition.

\newpage

\section{Point process of exceedances}
\label{sec:ppexceedances}

\subsection{Poisson approximation}

The main result of this section claims that the point process of exceedances converges to a Poisson point process in the transient case, i.e. $\alpha<1$. To introduce it, we denote for any $k\geq 1$ by \[\tau_k=\inf\{m\geq 0: \;\#\{S_1,\dots,S_m\}\geq k\}\] the time at which the random walk visits its $k$-th site. The \textit{point process of exceedances} is defined as 
\begin{equation}
\label{eq:defppexceedances}
  \Phi_n=\left\{\left(\frac{\tau_{k}}{n},\frac{\xi(S_{\tau_{k}})-b_{m(n)}}{ a_{m(n)}  }\right): \; \tau_k\leq n\right\}_{k\geq 1} \subset [0,1]\times \RR, 
\end{equation}
where $m(n)=\lfloor qn\rfloor$. 

\begin{Prop}{\label{processus1}}
	Let $\alpha < 1$. Assume that the sequence $(\xi(s))_{s\in \ZZ}$ satisfies the $\mathbf{D}(u_n)$ and $\mathbf{D'}(u_n)$ conditions for any threshold $u_n=u_{n}(x)=a_nx+b_n$, $x\in \RR$, satisfying Equation \eqref{eq:defnu}. Then $\Phi_n$ converges weakly to a Poisson point process $\Phi$ with intensity measure $m_{[0,1]}\otimes \nu$, where $m_{[0,1]}$ denotes the Lebesgue measure in $[0,1]$, i.e.
for any Borel subsets $B_1,\ldots, B_K \subset [0,1] \times \RR$ with $m_{[0,1]}\otimes \nu (\partial B_i)=0$, $1\leq i\leq K$, 
\[\left(\#\Phi_n\cap B_1, \ldots, \#\Phi_n\cap B_K \right) \overset{\mathcal{D}}{\conv[n]{\infty}}\left(\#\Phi\cap B_1, \ldots, \#\Phi\cap B_K \right).\] 
\end{Prop}

By using the Laplace functional, Franke and Saigo (Theorem 3 in \cite{franke_saigo_2009}) obtained a similar result when the $\xi(s)$'s are i.i.d. Proposition \ref{processus1} extends it and is based on Kallenberg's theorem. Our result is  stated only in the transient case, i.e. for $\alpha<1$. However, it remains true for $\alpha=1$ by taking $m(n) = \left\lfloor\frac{n}{h(n)}\right\rfloor$. When $\alpha>1$, the point process of exceedances is defined  in the same spirit as \eqref{eq:defppexceedances} by taking this time $m(n)=\lfloor n^{1/\alpha}\rfloor$. In this case, similarly to Theorem 4 in \cite{franke_saigo_2009}, we can show by adapting the proof of Proposition \ref{processus1} that $\Phi_n$ converges weakly to a Cox point process $\Phi_Y$, i.e. a Poisson point process in $[0,1]\times \RR$ with random intensity measure $\mu(\mathrm{d}t, \mathrm{d}x) = m_Y(\mathrm{d}t)\nu (\mathrm{d}x)$, where $m_Y(t)=m(Y(0,t))$.

\subsection{Technical results}
The proof of Proposition \ref{processus1} is mainly based on Kallenberg's theorem (see e.g. Proposition 3.22 in \cite{Resnick})  and on two technical lemmas which are stated below.

\begin{Th}[Kallenberg]
 Suppose $\Phi$ is a simple point process on $E$ and $\mathcal{I}$ is a basis of relatively compact open sets such that $\mathcal{I}$ is closed under finite unions and intersections and, for $I\in \mathcal{I}$, 
\[\PPP{\#\Phi \cap  \partial I=0}=1,\]
where $\partial I$ is the boundary of $I$.
Let $(\Phi_n)$ be a sequence of point processes on $E$ such that, for all $I\in \mathcal{I}$, 
\[\lim_{n\to +\infty} \mathbb{E}\left( \#\Phi_n \cap  I\right)=\mathbb{E}\left( \#\Phi \cap  I\right)\]
and 
\[\lim_{n\to +\infty}\PPP{\#\Phi_n \cap I=0}= \PPP{\#\Phi \cap I=0}.\]
Then $\Phi_n$ converges weakly to $\Phi$ in distribution.
\end{Th}

The following lemma is a direct adaptation of Lemma 1 in \cite{franke_saigo_2009} and  deals with the independence between the sequence $(\xi(S_n))_{n\geq 0}$ and the sequence $(\tau_k)_{k\geq 1}$.
\begin{Le}{\label{Idp_scenery_walk}}
	For all measurable sets $B\subset \mathbb{N}_+$ and $A\subset \mathbb{R}$, we have
\[\PPP{\tau_{k} \in B,\xi(S_{\tau_{k}})\in A}=\PPP{\tau_{k} \in B}\PPP{\xi\in A}.\]
\end{Le}

The second lemma is an extension of \cite{Nicolas_Ahmad}. More precisely,  under the assumptions that the $\mathbf{D}(u_{n})$ and $\mathbf{D'}(u_{n})$ conditions hold for the sequence $(\xi(s))_{s\in \ZZ}$,  we have shown in \cite{Nicolas_Ahmad} that 
\begin{equation*}
\PPP{\bigcap_{k \geq 1 : \frac{\tau_k}{n} \in (0,1]} \left\{\frac{\xi(S_{\tau_{k}})-b_{m(n)}}{a_{m(n)}}\notin (x,\infty)\right\}}  - \mathbb{E}\left(\exp\left(-\frac{R_n}{m(n)}\nu(x,\infty)\right)\right) \conv[n]{\infty}0
\end{equation*}
 when \eqref{eq:defnu} holds for any threshold $u_n=u_n(x)$, $x\in \RR$.  The following lemma deals with the case where the interval $(0,1]$ (resp. $(x,\infty))$ is replaced by $(a,b]$  (resp. $A\subset \RR)$ in the above equation.
 \begin{Le}{\label{th1_moi_nicolas}}
Let $A$ be a Borel subset in $\RR$ and let $0\leq a<b\leq 1$. Under the same assumptions as Proposition \ref{processus1}, for almost all realization of $(S_n)_{n\geq 0}$, we have
\begin{align*}
\lim\limits_{n\rightarrow\infty} \PPP{\bigcap_{k \geq 1 :\frac{\tau_{k}}{n} \in (a,b]} \left\{\frac{\xi(S_{\tau_{k}})-b_{m(n)}}{a_{m(n)}}\notin A\right\}} - \mathbb{E}\left(\exp\left(-\frac{R_{\lfloor nb \rfloor}-R_{\lfloor na\rfloor}}{m(n)}\nu(A)\right)\right)=0. 
\end{align*}
 \end{Le}

\subsection{Proofs}

\begin{prooft}{Lemma \ref{Idp_scenery_walk}}
Since the random walk and the random scenery are independent, we have 
\begin{align*}
\PPP{\tau_{k} \in B,\xi(S_{\tau_{k}})\in A}
&=\sum_{m \in B}\PPP{\tau_{k}=m,\xi(S_{m}) \in A }\notag \\
&=\sum_{m \in B}\sum_{s \in \mathbb{Z}}\PPP{\tau_{k}=m,S_{m}=s,\xi(s) \in A} \notag \\
&=\sum_{m \in B}\sum_{s \in \mathbb{Z}}\PPP{\tau_{k}=m,S_{m}=s}\PPP{\xi(s) \in A } \notag \\
&=\PPP{\tau_{k} \in B}\PPP{\xi \in A}.
\end{align*}
\end{prooft} 

 \begin{prooft}{Lemma \ref{th1_moi_nicolas}}
The proof will be sketched since it relies on a simple adaptation of the proof of Theorem 1 in \cite{Nicolas_Ahmad}.  

Let $(k_n)$, $(\ell_n)$ be as in  \eqref{eq:defrn} and let  
\begin{equation}
\label{eq:defrnbis}
r_n=\left\lfloor\frac{n}{k_n-1}\right\rfloor  + 1,
\end{equation} for $n$ large enough.  Given a realization of  $(S_n)_{n\geq 0}$, we write \[\mathcal{S}_{(na,nb]}=\left\{S_{\tau_k} : k\geq 1, \frac{\tau_k}{n} \in (a,b]\right\} \quad \text{and} \quad R_{\lfloor nb \rfloor}-R_{\lfloor na \rfloor}=\#\mathcal{S}_{(na,nb]}.\] To capture the fact that $(\xi(s))_{s\in \ZZ}$ satisfies the condition $\mathbf{D}(u_{n})$, we construct blocks and stripes as follows. 
Let \begin{equation*}  K_n=\left\lfloor\frac{R_{\lfloor nb \rfloor}-R_{\lfloor na \rfloor}}{r_n}\right\rfloor +1.\end{equation*}
We subdivide the set $\mathcal{S}_{(na,nb]}$ into subsets $B_i\subset \mathcal{S}_{(na,nb]}$, $1\leq i\leq K_n$, referred to as \textit{blocks}, in such a way that $\#B_i=r_n$ and $\max B_i<\min B_{i+1}$ for all $i\leq K_n-1$. Notice that  $K_n\leq k_n$ and $\#B_{K_n}=R_{\lfloor nb \rfloor}-R_{\lfloor na \rfloor}-(K_n-1)\cdot r_n$ a.s..  For each $j\leq K_n$, we denote by $L_j$ the family consisting of the $\ell_n$ largest  terms of $B_j$ (e.g. if $B_j=\{x_1,\ldots, x_{r_n}\}$, with $x_1<\cdots <x_{r_n}$, $j\leq K_n-1$, then $L_j=\{x_{r_n-\ell_n+1}, \ldots, x_{r_n}\}$). When $j=K_n$, we take the convention $L_{K_n}=\emptyset$ if $\#B_{K_n}<\ell_n$.  The set $L_{j}$  is referred to as a \textit{stripe}, and the union of the stripes is denoted by $\mathcal{L}_n=\bigcup_{j\leq K_n}L_{j}$. Proceeding in the same spirit as in the proofs of Lemmas 1 and 2 of \cite{Nicolas_Ahmad}, we can easily that for almost all realization of $(S_n)_{n\geq 0}$,
\begin{itemize}
\item $\PPP{\bigcap_{s \in \mathcal{S}_{(na,nb]}} \left\{\frac{\xi(s)-b_{m(n)}}{a_{m(n)}}\notin A\right\}}-\PPP{\bigcap_{s \in \mathcal{S}_{(na,nb]}\setminus \mathcal{L}_{n}} \left\{\frac{\xi(s)-b_{m(n)}}{a_{m(n)}}\notin A\right\}} \conv[n]{\infty} 0$;
\item $\PPP{\bigcap_{s \in \mathcal{S}_{(na,nb]}\setminus \mathcal{L}_{n}} \left\{\frac{\xi(s)-b_{m(n)}}{a_{m(n)}}\notin A\right\}}-\prod_{i\leq K_n}  \PPP{\bigcap_{s\in B_i\setminus \mathcal{L}_{n}} \left\{\frac{\xi(s)-b_{m(n)}}{a_{m(n)}}\notin A\right\}}\conv[n]{\infty} 0$;
    \item $\prod_{i\leq K_n}  \PPP{\bigcap_{s\in B_i\setminus \mathcal{L}_{n}} \left\{\frac{\xi(s)-b_{m(n)}}{a_{m(n)}}\notin A\right\}}-\prod_{i\leq K_n}  \PPP{\bigcap_{s\in B_i} \left\{\frac{\xi(s)-b_{m(n)}}{a_{m(n)}}\notin A\right\}} \conv[n]{\infty} 0$;
\item $\prod_{i\leq K_n}  \PPP{\bigcap_{s\in B_i} \left\{\frac{\xi(s)-b_{m(n)}}{a_{m(n)}}\notin A\right\}}-\mathbb{E}\left(\exp\left(-\frac{R_{\lfloor nb \rfloor}-R_{\lfloor na \rfloor}}{m(n)}\nu(A)\right)\right) \conv[n]{\infty} 0$. 
\end{itemize}
The first and the third assertions come from the fact that the size of the stripes is negligible compared to the size of the blocks, i.e. $\ell_n=o(r_n)$. The second assertion is a consequence of the fact that the sequence $(\xi(s))_{s \in \ZZ}$ satisfies the $\mathbf{D}(u_{n})$ condition and the last one is obtained by  using the $\mathbf{D}(u_{n})$ and $\mathbf{D'}(u_{n})$ conditions. Lemma \ref{th1_moi_nicolas} follows directly from the four assertions.
 \end{prooft}

\begin{prooft}{Proposition \ref{processus1}}
According to Kallenberg's theorem, it is sufficient to show that
\begin{enumerate}[(i)]
	\item $\lim\limits_{n\rightarrow\infty} \mathbb{E}\left(\#\Phi_n\cap I\right)=m_{[0,1]}\otimes\nu(I),$
	\item $\lim\limits_{n\rightarrow\infty} \PPP{\#\Phi_n\cap I=0}=e^{-m_{[0,1]}\otimes\nu(I)},$
\end{enumerate}
 for all set $I$ of the form $I=(a,b]\times A$, where $0\leq a<b\leq 1$ and where $A$ is an open subset of $E$. 

To deal with (i), we write
\begin{align*}
\mathbb{E}\left( \#\Phi_n \cap I\right) & =\sum_{k\geq 1}\PPP{\left(\frac{\tau_{k}}{n},\frac{\xi(S_{\tau_{k}})-b_{\lfloor qn \rfloor}}{a_{\lfloor qn \rfloor}}\right) \in I}\\
& = \sum_{k\geq 1}\PPP{\frac{\tau_{k}}{n} \in (a,b]}\PPP{\frac{\xi-b_{\lfloor qn \rfloor}}{a_{\lfloor qn \rfloor}} \in A} \notag \\
& = \sum_{k\geq 1} \PPP{\frac{\tau_{k}}{n} \in (a,b]}P_{\lfloor qn \rfloor}(A),
\end{align*} 
where the second line comes from  Lemma \ref{Idp_scenery_walk}. Using the fact that $\sum_{k\geq 1}\ind{\frac{\tau_k}{n}\in (a,b]}=R_{\lfloor nb \rfloor}-R_{\lfloor na \rfloor}$,  we have 
\begin{align*}
 \mathbb{E}\left(\#\Phi_n \cap I\right) &= \mathbb{E}\left(\sum_{k\geq 1} \ind{\frac{\tau_k}{n}\in (a,b]}\right) P_{\lfloor qn \rfloor}(A) \notag \\
&= \mathbb{E}\left(R_{\lfloor nb \rfloor}-R_{\lfloor na \rfloor}\right)P_{\lfloor qn \rfloor}(A).
\end{align*}
Moreover, according to Theorem \ref{th:legallrosen} and to the Lebesgue's dominated convergence theorem, we know that $\mathbb{E}(R_{\lfloor nb \rfloor}-R_{\lfloor na \rfloor})\eq[n]{\infty} nq(b-a)$.  This, together with \eqref{eq:defnu2} implies 
\begin{align*}
 \mathbb{E}\left(\#\Phi_n \cap I\right) \conv[n]{\infty} (b-a)\times \nu(A) = m_{[0,1]}\otimes \nu(I).  
\end{align*}

To deal with (ii), we observe that
\begin{align*}
 \PPP{\#\Phi_n\cap I=0}  =\PPP{\bigcap_{k \geq 1:\frac{\tau_{k}}{n} \in (a,b]} \left\{\frac{\xi(S_{\tau_{k}})-b_{\lfloor qn \rfloor}}{a_{\lfloor qn \rfloor}}\notin A\right\}}.
\end{align*}
According to Lemma \ref{th1_moi_nicolas}, Theorem \ref{th:legallrosen} and the Lebesgue's dominated convergence theorem, we have 
\begin{align*}
  \PPP{\#\Phi_n\cap I=0} & = \mathbb{E}\left(\exp\left(-\frac{R_{\lfloor nb\rfloor} - R_{\lfloor na\rfloor}}{\lfloor qn\rfloor}\nu(A)\right)\right) + o(1)\\
  & \conv[n]{\infty} \exp\left(-(b-a)\nu(A)\right).
  \end{align*}
This, together with the fact that $(b-a)\nu(A) = m_{[0,1]}\otimes \nu(I)$, concludes the proof of Proposition \ref{processus1}. 
\end{prooft}

\section{Properties of $(\xi(S_n))_{n\geq 0}$}
\label{sec:properties}
In this section, we give some properties of  $(\xi(S_n))_{n\geq 0}$. More precisely, we show that the latter satisfies the $D(u_n)$ condition and an extension of the so-called $D^{(k)}(u_n)$ condition, but does not satisfy the $D'(u_n)$ condition.

\subsection{Distributional mixing property}
The following extends Proposition 2 in \cite{franke_saigo_2009}, which deals with the case where the $\xi(s)$'s are i.i.d., to sequences which only satisfy the  $\mathbf{D}(u_n)$ and $\mathbf{D'}(u_n)$ conditions. 

\begin{Prop}
\label{th:xissatisfiesdun}
Let $\alpha < 1$. Assume that the sequence $(\xi(s))_{s\in \ZZ}$ satisfies the $\mathbf{D}(u_n)$ and $\mathbf{D'}(u_n)$ conditions for a threshold $u_n$ such that $n\PPP{\xi > u_n} \conv[n]{\infty}\tau$, with $\tau>0$. Then $(\xi(S_n))_{n\geq 0}$ satisfies the $\mathbf{D}(u_n)$ condition. 
\end{Prop}

\begin{prooft}{Proposition \ref{th:xissatisfiesdun}}
We adapt several arguments of \cite{franke_saigo_2009} in our context.  Let  
 $0\leq i_1< \dots <i_p<j_1<\dots <j_{p'} \leq n$ be a family of integers, with  $j_1-i_p>\ell_n$ and $k_n\ell_n=o(n)$. 
To prove that  $(\xi(S_n))_{n\geq 0}$ satisfies the $\mathbf{D}(u_n)$, we have to show that
  \[ |F'_{i_1,\dots,i_p,j_1,\dots,j_{p'}}(u_n) -F'_{i_1,\dots,i_p}(u_n)F'_{j_1,\dots,j_{p'}}(u_n)| \leq \tilde{\alpha}_{n,\ell_n},\]
  for some sequence $(\tilde{\alpha}_{n,\ell})_{(n,\ell)\in \NN^2}$ such that $k_n\tilde{\alpha}_{n,\ell_n}\conv[n]{\infty}0$,  with  \[F'_{i_1,\dots,i_p}(u_n) = \PPP{\xi(S_{i_1}) \leq u_n,\ldots, \xi(S_{i_p})\leq u_n}.\]
 We will use below the following notation: 
 \begin{itemize}
     \item $R_{i_1,\dots,i_p,j_1,\dots,j_{p'}}= \# \{S_{i_1},\dots ,S_{i_p},S_{j_1},\dots,S_{j_{p'}} \}$;
     \item $R_{i_1,\dots,i_p}=\# \{S_{i_1},\dots,S_{i_p}\}$;
     \item $ R_{j_1,\dots,j_{p'}}=\# \{S_{j_1},\dots,S_{j_{p'}}\}$;
     \item $R_{j_1,\dots,j_{p'}}^{i_1,\dots,i_p}=\# \{S_{i_1},\dots,S_{i_p}\} \cap \{S_{j_1},\dots,S_{j_{p'}}\}=R_{i_1,\dots,i_p}+R_{j_1,\dots,j_{p'}}-R_{i_1,\dots,i_p,j_1,\dots,j_{p'}}$. 
  \end{itemize}
  We have
 \begin{multline}
 \label{eq:majF'}
   |F'_{i_1,\dots,i_p,j_1,\dots,j_{p'}}(u_n) -F'_{i_1,\dots,i_p}(u_n)F'_{j_1,\dots,j_{p'}}(u_n)| \\
   \leq  \left\lvert F'_{i_1,\dots,i_p,j_1,\dots,j_{p'}}(u_n)-\mathbb{E}\left(\exp\left(-\frac{R_{i_1,\dots,i_p,j_1,\dots,j_{p'}}}{n}\tau\right)\right)\right\rvert \\
    + \left\lvert\mathbb{E}\left(\exp\left(-\frac{R_{i_1,\dots,i_p,j_1,\dots,j_{p'}}}{n}\tau\right)\right)-\mathbb{E}\left(\exp\left(-\frac{R_{i_1,\dots,i_p}+R_{j_1,\dots,j_{p'}}}{n}\tau\right)\right) \right\rvert \\
   +\left\lvert \mathbb{E}\left(\exp\left(-\frac{R_{i_1,\dots,i_p}+R_{j_1,\dots,j_{p'}}}{n}\tau\right)\right)-F'_{i_1,\dots,i_p}(u_n)F'_{j_1,\dots,j_{p'}}(u_n)\right\rvert. 
 \end{multline}
To deal with the first and the third terms of the right-hand side of \eqref{eq:majF'}, we will use the following lemma.

\begin{Le}
\label{le:deviationrange}
For almost all realization of $(S_n)_{n\geq 0}$ and for all $0\leq i_1<i_2<\cdots <i_p\leq n$,
\[\left| F'_{i_1,\ldots, i_p}(u_n)   -  \exp\left( - \frac{R_{i_1,\ldots, i_p}}{n}\tau    \right)  \right|   \leq \varepsilon_n,\]
with $\varepsilon_n = \varepsilon_n^{(1)}+\varepsilon_n^{(2)}$, where  $\varepsilon_n^{(1)}$ and $\varepsilon_n^{(2)}$ are defined in \eqref{eq:defepsilon1} and \eqref{eq:defepsilon2} respectively. 
\end{Le}

\begin{prooft}{Lemma \ref{le:deviationrange}}
Similarly to Lemma \ref{th1_moi_nicolas}, the main idea is to adapt several arguments appearing in the proofs of Lemmas 1 and 2 in \cite{Nicolas_Ahmad} in our context. Let $(k_n)$ and $(r_n)$ be as in  \eqref{eq:defrn} and \eqref{eq:defrnbis}. Given $1\leq i_1<i_2<\cdots <i_p\leq n$, we subdivide the random set $\{S_{i_1},\ldots, S_{i_p}\}$ into $K_n$ blocks, with $K_n=\lfloor \frac{R_{i_1,\ldots, i_p}}{r_n} \rfloor + 1$, in the same spirit as we did in the proof of Lemma \ref{th1_moi_nicolas}.  More precisely, there exists a unique $K_n$-tuple of subsets $B_i\subset \mathcal{S}_n$, $i\leq K_n$, such that the following properties hold: $\bigcup_{j\leq K_n}B_j= \{S_{i_1},\ldots, S_{i_p}\}  $, $\#B_i=r_n$ and $\max B_i<\min B_{i+1}$ for all $i\leq K_n-1$. In particular, we have  $K_n\leq k_n$ and $\#B_{K_n}=R_n-(K_n-1)\cdot r_n$ a.s.. Without loss of generality, we assume that $\#B_{K_n} = \#B_i = r_n$ for all $i\leq K_n-1$, so that $R_{i_1,\ldots, i_p} = K_nr_n$. For each $j\leq K_n$, we also denote by $L_j$ the family consisting of the $\ell_n$ largest  terms of $B_j$ and we let $\mathcal{L}_n=\bigcup_{j\leq K_n}L_{j}$. In the rest of the paper, we write $M_B=\max_{s\in B}\xi(s)$ for all subset $B\subset \ZZ$. 

Adapting the proof of Lemma 1 in \cite{Nicolas_Ahmad}, we can show that the following inequalities hold for almost all realization of $(S_n)_{n\in \geq 0}$ and for $n$ larger than some deterministic integer $n_0$:
\begin{equation*}
\left|  \PPP{M_{\{S_{i_1}, \ldots, S_{i_p}\}}  \leq u_n  } -  \PPP{M_{\{S_{i_1}, \ldots, S_{i_p}\} \setminus \mathcal{L}_n }  \leq u_n  }    \right| \leq k_n \ell_n \PPP{\xi>u_n};
\end{equation*}
\begin{equation*}
\left|  \PPP{M_{\{S_{i_1}, \ldots, S_{i_p}\} \setminus \mathcal{L}_n }  \leq u_n  } - \prod_{j\leq K_n}\PPP{M_{B_j\setminus \mathcal{L}_n} \leq u_n}   \right| \leq k_n \alpha_{n, \ell_n};
\end{equation*}
\begin{equation*}
\left|  \prod_{j\leq K_n}\PPP{M_{B_j\setminus \mathcal{L}_n} \leq u_n}  -  \prod_{j\leq K_n}\PPP{M_{B_j}\leq u_n}    \right| \leq 2\frac{\tau k_n\ell_n}{n}.
\end{equation*}
Since $ F'_{i_1,\ldots, i_p}(u_n) =  \PPP{M_{\{S_{i_1}, \ldots, S_{i_p}\}}  \leq u_n  } $ and $\PPP{\xi>u_n}\eq[n]{\infty}\frac{\tau}{n}$,  we get for almost all realization of $(S_n)_{n\geq  0}$, 
\[\left|   F'_{i_1,\ldots, i_p}(u_n) -  \prod_{j\leq K_n}\PPP{M_{B_j}\leq u_n}    \right|\leq \varepsilon_n^{(1)},\]  with
 \begin{equation}
\label{eq:defepsilon1} 
 \varepsilon_n^{(1)}  = c\cdot \left(\frac{k_n\ell_n}{n} + k_n \alpha_{n,\ell_n}\right).
 \end{equation}

Without loss of generality, we assume from now on that $\PPP{\xi>u_n}=\frac{\tau}{n}$. We show below that
\begin{equation}
\label{eq:aimrange}
\left|   \prod_{j\leq K_n}\PPP{M_{B_j}\leq u_n}  -  \exp\left( - \frac{R_{i_1,\ldots, i_p}}{n}\tau    \right)   \right|\leq \varepsilon_n^{(2)},
\end{equation} for some  deterministic sequence $\varepsilon_n^{(2)}\conv[n]{\infty}0$. To do it, we adapt several arguments of Lemma 2 in \cite{Nicolas_Ahmad}. First, we notice that for $n$ large enough, 
\begin{multline*}
 \prod_{j\leq K_n}\PPP{M_{B_j}\leq u_n} -  \exp\left( - \frac{R_{i_1,\ldots, i_p}}{n}\tau    \right) \\
 \begin{split}
 & \geq   \exp\left( K_n\log (1-r_n\PPP{\xi>u_n})   \right) -  \exp\left( - \frac{R_{i_1,\ldots, i_p}}{n}\tau    \right)\\
 & \geq \exp\left(  -K_nr_n\PPP{\xi>u_n} - K_n(r_n\PPP{\xi >u_n})^2  \right) - \exp\left( - \frac{R_{i_1,\ldots, i_p}}{n}\tau    \right),
 \end{split}
\end{multline*}
where the last line comes from the facts that $\log(1-x)\geq -x-x^2$ for $|x|$ small enough and that $r_n\PPP{\xi>u_n}\conv[n]{\infty}0$. Because $K_nr_n=R_{i_1,\ldots, i_p}$ and $\PPP{\xi>u_n}=\frac{\tau}{n}$, we have
\begin{multline*}
 \prod_{j\leq K_n}\PPP{M_{B_j}\leq u_n} -  \exp\left( - \frac{R_{i_1,\ldots, i_p}}{n}\tau    \right) \\
 \begin{split}
 & \geq  \exp\left( - \frac{R_{i_1,\ldots, i_p}}{n}\tau    \right) \left( \exp\left(  -K_n (r_n\PPP{\xi >u_n})^2 \right)  - 1  \right)\\
 & \geq \exp(-k_n (r_n\PPP{\xi >u_n})^2) - 1,
 \end{split}
\end{multline*}
where the last line comes from the fact that $K_n\leq k_n$ a.s.. Since $k_nr_n \eq[n]{\infty}n$, we have
\[ \prod_{j\leq K_n}\PPP{M_{B_j}\leq u_n} -  \exp\left( - \frac{R_{i_1,\ldots, i_p}}{n}\tau    \right) \geq c\cdot \frac{1}{k_n}.\]

Moreover, because $ \prod_{j\leq K_n}\PPP{M_{B_j}\leq u_n} \leq \exp\left( - \sum_{j\leq K_n} \PPP{M_{B_j}>u_n} \right)$, it follows from the Bonferroni inequalities (see e.g. p110 in Feller \cite{Fel}) that 
\begin{multline*}
 \prod_{j\leq K_n}\PPP{M_{B_j}\leq u_n}\\
\leq \exp\left( -{(K_n-1)}r_n \PPP{\xi>u_{n}}{+}\sum_{j\leq K_n} \sum_{\alpha < \beta; \alpha, \beta \in B_{j}}\PPP{\xi(\alpha)>u_{n},\xi(\beta)>u_{n}} \right).
\end{multline*}
Since $K_nr_n=R_{i_1,\ldots, i_p}$ and $\PPP{\xi>u_n} = \frac{\tau}{n}$, we have
\begin{multline*}
 \prod_{j\leq K_n}\PPP{M_{B_j}\leq u_n}-  \exp\left( - \frac{R_{i_1,\ldots, i_p}}{n}\tau    \right) = \exp\left( - \frac{R_{i_1,\ldots, i_p}}{n}\tau    \right)\\
 \times  \left( \exp\left( r_n\PPP{\xi >u_n} +  \sum_{j\leq K_n} \sum_{\alpha < \beta; \alpha, \beta \in B_{j}}\PPP{\xi(\alpha)>u_{n},\xi(\beta)>u_{n}}   \right) - 1\right)
 \end{multline*}
 and therefore
\begin{multline*}
 \prod_{j\leq K_n}\PPP{M_{B_j}\leq u_n}-  \exp\left( - \frac{R_{i_1,\ldots, i_p}}{n}\tau    \right) \\
  \leq \exp\left( r_n\PPP{\xi >u_n} +  \sum_{j\leq K_n} \sum_{\alpha < \beta; \alpha, \beta \in B_{j}}\PPP{\xi(\alpha)>u_{n},\xi(\beta)>u_{n}}   \right) - 1.
 \end{multline*}  
Proceeding along the same lines as in the proof of Lemma 2 in \cite{Nicolas_Ahmad}, we can show that 
\begin{multline*}
\exp\left( r_n\PPP{\xi >u_n} +  \sum_{j\leq K_n} \sum_{\alpha < \beta; \alpha, \beta \in B_{j}}\PPP{\xi(\alpha)>u_{n},\xi(\beta)>u_{n}}   \right) - 1\\
\leq c\left( \frac{1}{k_n} + n \sum_{s=1}^{\lfloor n/k_n\rfloor} \PPP{\xi(0)>u_n, \xi(s)>u_n} \right).
\end{multline*}
This shows \eqref{eq:aimrange} with 
\begin{equation}
\label{eq:defepsilon2} 
\varepsilon_n^{(2)} = c\left( \frac{1}{k_n} + n \sum_{s=1}^{\lfloor n/k_n\rfloor} \PPP{\xi(0)>u_n, \xi(s)>u_n} \right).
\end{equation}
and consequently concludes the proof of Lemma \ref{le:deviationrange}. 
\end{prooft} 
According to \eqref{eq:defrn}, the fact that $(\xi(s))_{s\in \ZZ}$ satisfies the $\mathbf{D}'(u_n)$ condition and the fact that  $k_n\alpha_{n,\ell_n}\conv[n]{\infty}0$, we have $\varepsilon_n\conv[n]{\infty}0$. It follows from Lemma \ref{le:deviationrange} that the first and the third terms of the right-hand side of \eqref{eq:majF'} converge to 0 as $n$ goes to infinity. To deal with the second one, we write
 \begin{multline*}
   \left|  \exp\left(-\frac{R_{i_1,\dots,i_p,j_1,\dots,j_{p'}}}{n}\tau\right)-\exp\left(-\frac{R_{i_1,\dots,i_p}+R_{j_1,\dots,j_{p'}}}{n}\tau\right) \right| \\
   \begin{split}
     &=\exp\left(-\frac{R_{i_1,\dots,i_p}+R_{j_1,\dots,j_{p'}}}{n}\tau\right) \left(\exp\left(\frac{R_{i_1,\dots,i_p}^{j_1,\dots,j_{p'}}}{n}\tau\right)-1\right) \\
     & \leq \exp\left(\frac{R_{1,\dots,i_p}^{i_p+\ell_n+1,\ldots, n}}{n}\tau\right)-1,
     \end{split}
 \end{multline*}
 where the last line comes from the fact that $j_{1}-i_p > \ell_n$. Since $\ell_n\geq 0$, we get 
\begin{multline}
 \label{eq:majrangeinter}
 \sup \left|  \exp\left(-\frac{R_{i_1,\dots,i_p,j_1,\dots,j_{p'}}}{n}\tau\right)-\exp\left(-\frac{R_{i_1,\dots,i_p}+R_{j_1,\dots,j_{p'}}}{n}\tau\right) \right|\\
 \leq \sup_{i\leq n}  \exp\left(\frac{R_{1,\dots,i}^{i+1,\ldots, n}}{n}\tau\right)-1 ,
\end{multline} 
 where the supremum in the left-hand side is taken over all integers $0\leq i_1< \dots <i_p<j_1<\dots <j_{p'} \leq n$, with  $j_1-i_p>\ell_n$. Moreover, using the fact that $R_{1,\dots,i}^{i+1,\ldots, n} = R_{1,\ldots, i}+R_{i+1,\ldots, n}  -R_{1,\ldots, n}$ and following \cite{LeGRo}, we have $\sup_{i\leq n}\frac{R_{1,\dots,i}^{i+1,\ldots, n}}{n}\conv[n]{\infty}0$ a.s..  This, together with  \eqref{eq:majrangeinter} and the Lebesgue's dominated convergence theorem implies
\[\sup \left| \EEE{ \exp\left(-\frac{R_{i_1,\dots,i_p,j_1,\dots,j_{p'}}}{n}\tau\right)} - \EEE{\exp\left(-\frac{R_{i_1,\dots,i_p}+R_{j_1,\dots,j_{p'}}}{n}\tau\right)} \right|\conv[n]{\infty}0\]
and consequently concludes the proof of Proposition \ref{th:xissatisfiesdun}. 
\end{prooft}


\subsection{The $D^{(k)}(u_n)$ as $k\rightarrow\infty$}

 In \cite{chernick_hsing_mccormick_1991}, the authors introduce a local mixing condition, referred to as the $D^{(k)}(u_n)$ condition, which allows to express the extremal index in terms of joint distribution. We recall the latter below.

\paragraph{Condition $\mathbf{D}^{(k)}(u_n)$} Let $(\xi(s))_{s \in \ZZ}$ be a sequence of random variables and let $u_n$ be a threshold such that $n\PPP{\xi>u_n}\conv[n]{\infty}\tau$, for some $\tau>0$. In conjunction with the $D(u_n)$ condition, we say that the $D^{(k)}(u_n)$ condition, $k\geq 1$, holds if there exist two sequences of integers $(k_n)$ and $(\ell_n)$ such that
\begin{equation*}
   k_n \rightarrow \infty, \; \quad k_n \alpha_{n,\ell_n}\rightarrow 0, \;  k_n\ell_n=o(n) 
\end{equation*}
and 
\begin{equation}{\label{conditionD^k}}
  \lim\limits_{n\rightarrow\infty} n \PPP{\xi({1})>u_n \geq M_{2,k}, \; M_{k+1,r_n}>u_n}=0,
\end{equation}
where $r_n$ is as in \eqref{eq:defrnbis}  and where $M_{i,j} = \max\{\xi(i), \xi(i+1),\ldots, \xi(j)\}$ for all $i\leq j$, with the convention $M_{i,j}=-\infty$ if $i>j$. 
As mentioned in \cite{chernick_hsing_mccormick_1991}, Equation \eqref{conditionD^k} is implied by the condition 
\begin{align*}
  \lim\limits_{n\rightarrow\infty} n \sum_{s=k+1}^{r_n} \PPP{\xi({1})>u_n \geq M_{2,k}, \; \xi(s)>u_n}=0.
\end{align*}
Observe that the last line is the $D'(u_n)$ condition if $k=1$. 

Roughly, the following proposition states that the sequence $(\xi(S_n))_{n\geq 0}$ satisfies the $D^{(k)}(u_n)$ condition as $k$ goes to infinity.
\begin{Prop}{\label{D_infini}}
Under the same assumptions as Proposition \ref{th:xissatisfiesdun}, we have
\[\lim\limits_{k\rightarrow\infty} \lim\limits_{n\rightarrow\infty}  n \sum_{j=k+1}^{r_n} \PPP{\xi({S_1})>u_n \geq M'_{2,k}, \; \xi(S_j)>u_n}=0,  
\]
where $M'_{i,j}= \max_{i\leq t \leq j} \xi(S_{t})$ if  $ i\leq j$ and  $M'_{i,j} = -\infty$ if  $i>j$.
\end{Prop}
\begin{prooft}{Proposition \ref{D_infini}} 
 For all  $k\geq 1$, we have 
 \begin{multline}
 \label{eq_D_D'_D^infini}
     n \sum_{j=k+1}^{r_n} \PPP{\xi(S_{1})>u_n \geq M'_{2,k}, \; \xi(S_j)>u_n} \\
     =n \sum_{j=k+1}^{r_n} \PPP{\xi(S_{1})>u_n \geq M'_{2,k}, \; \xi(S_j)>u_n|S_j=S_1}\PPP{S_j=S_1}  \\
     +n \sum_{j=k+1}^{r_n} \PPP{\xi(S_{1})>u_n \geq M'_{2,k}, \; \xi(S_j)>u_n|S_j\neq S_1}\PPP{S_j\neq S_1}.
 \end{multline}
 The first term of the right-hand side of \eqref{eq_D_D'_D^infini}  tends to zero as $k,n\rightarrow \infty$. Indeed,  
 \begin{align*}
      \PPP{\xi(S_{1})>u_n \geq M'_{2,k}, \; \xi(S_j)>u_n|S_j=S_1}&\PPP{S_j=S_1} \\
     &\leq   \PPP{\xi(S_{1})>u_n }\PPP{S_j=S_1}.
 \end{align*}
Moreover, because $(S_n)_{n\geq 1}$ is a transient random walk, we have  $\sum_{j=2}^{\infty}\PPP{S_j=S_1}< \infty$, which implies \[\lim\limits_{k\rightarrow\infty} \lim\limits_{n\rightarrow\infty}n \PPP{\xi(S_{1})>u_n}\sum_{j=k+1}^{r_n}\PPP{S_j=S_1}=0,\]
and therefore 
\[\lim\limits_{k\rightarrow\infty} \lim\limits_{n\rightarrow\infty} n \sum_{j=k+1}^{r_n} \PPP{\xi(S_{1})>u_n \geq M'_{2,k}, \; \xi(S_j)>u_n|S_j=S_1}\PPP{S_j=S_1} =0.\]

To prove that the second term of the right-hand side of \eqref{eq_D_D'_D^infini} goes to $0$, we write 
 \begin{multline}{\label{boule}}
  n \sum_{j=k+1}^{r_n} \PPP{\xi(S_{1})>u_n \geq M'_{2,k},\  \xi(S_j)>u_n|S_j\neq S_1}\PPP{S_j\neq S_1}  \\
  = n \sum_{j=k+1}^{r_n} \PPP{\xi(S_{1})>u_n \geq M'_{2,k},\  \xi(S_j)>u_n|S_j \in B^*(S_1, r_n)} \PPP{S_j \in B^*(S_1, r_n)} \\
  + n \sum_{j=k+1}^{r_n} \PPP{\xi(S_{1})>u_n \geq M'_{2,k},\  \xi(S_j)>u_n|S_j \notin B(S_1, r_n)} \PPP{S_j \notin B(S_1, r_n)},
 \end{multline}
 where $B(S_1, r_n):=\left\{S \in \mathcal{S}_n: \, |S-S_1|\leq r_n \right\}$ and $B^*(S_1, r_n)=B(S_1, r_n) \setminus \{S_1\}$. We prove below that the last two terms in \eqref{boule} converge to $0$. For the first one, we write
 \begin{multline*}
  n \sum_{j=k+1}^{r_n} \PPP{\xi(S_{1})>u_n \geq M'_{2,k},\  \xi(S_j)>u_n|S_j \in B^*(S_1, r_n)} \PPP{S_j \in B^*(S_1, r_n)}\\
   \leq n \sum_{j=2}^{r_n} \PPP{\xi(0)>u_n ,\  \xi(S_j-S_1)>u_n|S_j \in B^*(S_1, r_n)}.
 \end{multline*}
 The last quantity converges to $0$ as $n$ goes to infinity since the sequence $(\xi(s))_{s \in \ZZ}$ satisfies the $\mathbf{D'}(u_{n})$ condition. To deal with the second term of \eqref{boule}, we write
 \begin{multline*}{\label{outsideBoule}}
  n\sum_{j=k+1}^{r_n} \PPP{\xi(S_{1})>u_n \geq M'_{2,k},\  \xi(S_j)>u_n|S_j \notin B(S_1, r_n)} \PPP{S_j \notin B(S_1, r_n)} \\
  \begin{split}
  &  \leq n \sum_{j=k+1}^{r_n} \PPP{\xi(S_{1})>u_n,\  \xi(S_j)>u_n|S_j \notin B(S_1, r_n)}  \\ 
 & \leq n\sum_{j=k+1}^{r_n}\PPP{\xi>u_n}^2+n\sum_{j=k+1}^{r_n}\left|\PPP{\xi(S_{1})>u_n, \;  \xi(S_j)>u_n|S_j \notin B(S_1, r_n)}- \PPP{\xi>u_n}^2\right|.
 \end{split}
 \end{multline*}
The first series tends to 0 as $n$ goes to infinity because 
 \[n\sum_{j=k+1}^{r_n}\PPP{\xi>u_n}^2\leq nr_n\PPP{\xi>u_n}^2\eq[n]{\infty}\tau^2\frac{r_n}{n},\]
and  $r_n=o(n)$. To deal with the second series, we use the $\mathbf{D}(u_{n})$ condition. This gives 
 \begin{align*}
     n\sum_{j=k+1}^{r_n} |\PPP{\xi(S_{1})>u_n, \xi(S_j)>u_n|S_j \notin B(S_1, r_n)}- \PPP{\xi>u_n}^2|   &\leq nr_n \alpha_{n,r_n}  \\
     &\leq \frac{n^2}{k_n}\alpha_{n,r_n},
 \end{align*}
 which converges to $0$ as $n$ goes to infinity according to \eqref{eq:defrn}. This concludes the proof of Proposition \ref{D_infini}.
\end{prooft}

\subsection{The extremal index}
\label{sec:extremalindex}

Let $(k_n)$ and $(r_n)$ be as in \eqref{eq:defrn} and \eqref{eq:defrnbis}. Let us denote by $R_n=\#\mathcal{S}_n$ and $K_n = \left\lfloor \frac{R_n}{r_n} \right\rfloor + 1$. The following proposition deals with $M_{\mathcal{S}_n}$ under the $\mathbf{D}(u_n)$ condition.

\begin{Prop}{\label{ob_moi}}
Let $\alpha < 1$. Assume that the sequence $(\xi(s))_{s\in \ZZ}$ satisfies the $\mathbf{D}(u_n)$ conditions for a threshold $u_n$ such that $n\PPP{\xi > u_n} \conv[n]{\infty}\tau$, with $\tau>0$. Then for almost all realization of  $(S_n)_{n \geq 0}$, 
\[\PPP{M_{\mathcal{S}_n}\leq u_n}-\exp\left(-\sum_{j=1}^{K_n} \sum_{i=1}^{r_n}\PPP{\xi(S_{((j-1)r_n+i)})>u_n\geq M'_{((j-1)r_n+i+1, \; jr_n)}}\right) \conv[n]{\infty}0,\]
where \[M'_{(i,j)}:= \left\{ \begin{split} & \max_{i\leq t \leq j} \xi(S_{(t)}), \quad  i\leq j \\
& -\infty, \quad  i>j \end{split} \right.\] and where $S_{(t)}$ is the $t$-th largest value of the $\xi(S_i)'s$, $i\leq n$.
\end{Prop}

A similar result was obtained by O'Brien (Theorem 2.1. in  \cite{Ob}). However, the above proposition is not a consequence of the latter. Proposition \ref{ob_moi} remains true if the sequence $(\xi(s))_{s\in \ZZ}$ only satisfies the $D(u_n)$ condition (i.e. when $k_n\alpha_{n,\ell_n}\conv[n]{\infty}0$ instead of $\frac{n^2}{k_n}\alpha_{n, \ell_n}\conv[n]{\infty}0$). As a direct consequence of such a result, if for almost all realization of $(S_n)_{n\geq 0}$, 
   \begin{equation*}
   \label{eq:condtheta}
   \frac{ 1}{n}\sum_{j \leq K_n} \sum_{i=1}^{r_n}\PPP{M'_{((j-1)r_n+i+1, \; jr_n)}\leq u_n| \xi(S_{((j-1)r_n+i)})>u_n} \conv[n]{\infty} \theta,
   \end{equation*} for some $\theta \in [0,1]$, then $\PPP{M_{\mathcal{S}_n} \leq u_n} \conv[n]{\infty} e^{-\theta\tau}$. In this case, the term $\theta$ is referred to as the \textit{extremal index} (see e.g.  \cite{L2}) and can be interpreted as the reciprocal of the mean size of a cluster of exceedances. As stated in Theorem 1 in \cite{Nicolas_Ahmad}, when the sequence $(\xi(s))_{s\in \ZZ}$ satisfies the $\mathbf{D}(u_{n})$ and $\mathbf{D'}(u_{n})$ conditions, we have
\begin{equation}
\label{eq:extremalindexq}
\PPP{M_{\mathcal{S}_n}\leq u_n}\conv[n]{\infty}e^{-q\tau }.
\end{equation}   
 In other words, under these conditions, the extremal index $\theta$ exists and  $\theta=q$.

\begin{prooft}{Proposition \ref{ob_moi}}
Let us write $\mathcal{S}_n=\{S_{(1)}, \dots,S_{(R_n)}\}$ with $S_{(1)}<S_{(2)}<\dots<S_{(R_n)}$, and partitition $\mathcal{S}_n$ into $K_n$ blocks  as in Lemma \ref{th1_moi_nicolas}. Without loss of generality, assume that the last block has the same size as the others, so that $\frac{R_n}{K_n}$ is an integer. Let $B_j=\{S_{((j-1)r_n+1)},\dots,S_{(jr_n)}\}$ be the $j$-th block of size $r_n$.  According to Lemma 1 in \cite{Nicolas_Ahmad}, for almost all realization of $(S_n)_{n\geq 0}$, we have  
\begin{equation*}
  \PPP{M_{\mathcal{S}_n}\leq u_n} -\exp\left(\sum_{j\leq K_n}\log \left(1-\PPP{M_{B_j}>u_n}\right)\right)\conv[n]{\infty} 0.
\end{equation*}
Moreover, because $| \log(1-x)+x | \leq Cx^2$ for $|x|$ small enough and  because $\PPP{M_{B_j}>u_n}\leq r_n \PPP{\xi>u_n}$  converges to $0$ as $n$ goes to infinity, we have 
\begin{multline*}
  \left\lvert\sum_{j\leq K_n}\log \left(1-\PPP{M_{B_j}>u_n}\right)+\sum_{j\leq K_n} \PPP{M_{B_j}>u_n} \right\rvert  \\
\begin{split} 
& \leq \sum_{j\leq K_n}  \left\lvert\log \left(1- \PPP{M_{B_j}>u_n}\right)+ \PPP{M_{B_j}>u_n} \right\rvert \\
& \leq C \sum_{j\leq K_n} \PPP{M_{B_j}>u_n}^2 \\
& \leq C k_n r_n^2 \PPP{\xi>u_n}^2.
\end{split}
\end{multline*}
The last term converges to $0$ as $n$ goes to infinity since $k_nr_n \eq[n]{\infty} n$, $n\PPP{\xi>u_n} \conv[n]{\infty} \tau$ and $r_n\PPP{\xi>u_n} \conv[n]{\infty} 0$. This shows that for almost all realization of $(S_n)_{n\geq 0}$ 
\begin{equation}
\label{eq:obrien2} 
\PPP{M_{\mathcal{S}_n}\leq u_n} -\exp\left(-\sum_{j\leq K_n} \PPP{M_{B_j}>u_n}\right)\conv[n]{\infty} 0.
\end{equation}
Besides, following the same lines as \cite{Ob}, we have
\begin{align*}
  \PPP{M_{B_j} \leq u_n} &= 1- \PPP{M_{B_j}>u_n}\notag \\
  & = 1- \sum_{i=1}^{r_n} \PPP{\xi(S_{((j-1)r_n+i)})>u_n\geq M'_{((j-1)r_n+i+1, \; jr_n)}} 
\end{align*}
This together with \eqref{eq:obrien2} concludes the proof of Proposition \ref{ob_moi}.
\end{prooft}

\subsection{The $D'(u_n)$ condition}
Recall that, in the classical literature  (see e.g. (3.2.1) in \cite{lucarini}), the $D'(u_n)$ condition holds for the sequence $(Z_n)$ if, in conjunction with the $D(u_n)$ condition,
\[\lim\limits_{n\rightarrow \infty} n \sum_{i=2}^{[ n/k_n]}\PPP{Z_1>u_n, Z_i>u_n}=0,\]
for some sequence of integers $(k_n)$ such that $k_n\conv[n]{\infty}\infty$, $k_n\alpha_{n,\ell_n}\conv[n]{\infty}0$ and $k_n\ell_n=o(n)$.  The following result is an extension of Proposition 3 in \cite{franke_saigo_2009}. However, we give  a simpler proof which is based on \cite{L2}.

\begin{Prop}
\label{prop:d'notsatisfied}
Under the same assumptions as Proposition \ref{th:xissatisfiesdun}, the sequence $(\xi(S_n))_{n\geq 0}$ does not satisfy the $D'(u_n)$ condition. 
\end{Prop}

\begin{prooft}{Proposition \ref{prop:d'notsatisfied}}
On the opposite, if $(\xi(S_n))_{n\geq 0}$ satisfies the $D'(u_n)$ condition, then $\PPP{M_{\mathcal{S}_n}\leq u_n}\conv[n]{\infty}e^{-\tau}$ according to Theorem 1.2 in \cite{L2}. This contradicts \eqref{eq:extremalindexq} since $q\neq 1$. 
\end{prooft}



\end{document}